\documentclass{gtpart}     
\usepackage{graphicx}  
\usepackage{amsmath}
\usepackage{amssymb}
\usepackage{verbatim}
\usepackage{amsfonts}
\title{Small curvature laminations in hyperbolic 3-manifolds}

\author{William Breslin}
\givenname{William}
\surname{Breslin}
\address{Department of Mathematics, University of Michigan, Ann Arbor, MI 48109}
\email{breslin@umich.edu}
\urladdr{http://www-personal.umich.edu/~breslin/index.html}


\arxivreference{0901.1330}
\arxivpassword{zddge}

\volumenumber{}
\issuenumber{}
\publicationyear{}
\papernumber{}
\startpage{}
\endpage{}
\doi{}
\MR{}
\Zbl{}
\received{}
\revised{}
\accepted{}
\published{}
\publishedonline{}
\proposed{}
\seconded{}
\corresponding{}
\editor{}
\version{}

\newtheorem{thm}{Theorem}
\theoremstyle{definition}
\newcommand{\h}{\mathbb{H}}
\newcommand{\bd}{\partial}

\begin{document}

\begin{abstract}    
We show that if $\mathcal{L}$ is a codimension-one lamination in a finite volume hyperbolic 3-manifold  such that the principal curvatures of each leaf of $\mathcal{L}$ are all in the interval $(-\delta ,\delta)$ for a fixed $\delta\in[0,1)$ and no complementary region of $\mathcal{L}$ is an interval bundle over a surface, then each boundary leaf of $\mathcal{L}$ has a nontrivial fundamental group.  We also prove existence of a fixed constant $\delta_0 > 0$ such that if $\mathcal{L}$ is a codimension-one lamination in a finite volume hyperbolic 3-manifold such that the principal curvatures of each leaf of $\mathcal{L}$ are all in the interval $(-\delta_0 ,\delta_0)$ and no complementary region of $\mathcal{L}$ is an interval bundle over a surface, then each boundary leaf of $\mathcal{L}$ has a noncyclic fundamental group.
\end{abstract}

\maketitle


\section{Introduction}

In \cite{Zeghib}, Zeghib proved that any totally geodesic codimension-one lamination in a closed hyperbolic 3-manifold is a finite union of disjoint closed surfaces. In this paper we investigate whether a similar result holds for codimension-one laminations with small principal curvatures. We will prove the following theorems:

\begin{thm}\label{thm1}
Let $\delta \in [0,1)$.  If $\mathcal{L}$ is a codimension-one lamination in a finite volume hyperbolic 3-manifold such that the principal curvatures of each leaf of $\mathcal{L}$ are everywhere in $(-\delta ,\delta )$ for a fixed constant $\delta\in[0,1)$ and no complementary region of $\mathcal{L}$ is an interval bundle over a surface, then each boundary leaf of $\mathcal{L}$ has a nontrivial fundamental group.
\end{thm}

\begin{thm}\label{thm2}
There exists a fixed constant $\delta_0 > 0$ such that if $\mathcal{L}$ is a codimension-one lamination in a finite volume hyperbolic 3-manifold such that the principal curvatures of each leaf of $\mathcal{L}$ are everywhere in $(-\delta_0 ,\delta_0 )$ and no complementary region is an interval bundle over a surface, then each boundary leaf of $\mathcal{L}$ has a noncyclic fundamental group.
\end{thm}

\section{Examples}

Let $\mathcal{L}$ be a codimension-one lamination in a complete hyperbolic 3-manifold $M$.  Let $L$ be a leaf of $\mathcal{L}$ and endow it with the path metric induced from $M$.  Let $\tilde{L}$ be the universal cover of $L$ and lift the inclusion $i_L: L \rightarrow M$ to a map $\tilde{i_L} : \tilde{L} \rightarrow \h^3$.  A map $f: X \rightarrow Y$ from a metric space $X$ to a metric space $Y$ is a $(k,c)$-quasi-isometry if $\frac{1}{k} d_{X}(a,b) - c \le d_{Y}(f(a),f(b)) \le k d_{X}(a,b) + c$.   The leaf $L$ is \textit{quasi-isometric} if $\tilde{i}_L$ is a $(k,c)$-quasi-isometry for some $k,c$.  The lamination $\mathcal{L}$ is \textit{quasi-isometric} if each leaf of $\mathcal{L}$ is quasi-isometric for the same fixed constants $k,c$.

Let $\delta \in (0,1)$. If the principal curvatures of $\tilde{i}_L (\tilde{L})$ are everywhere in $(-\delta ,\delta)$, then the map $\tilde{i}_L$ is a $(k,c)$-quasi-isometry for constants $k,c$ depending only on $\delta$ (see Thurston \cite{tnotes}.  Also see Leininger \cite{Leininger} for an elementary proof).

The constant $\delta_0$ in Theorem \ref{thm2} is less than 1, so a lamination satisfying the hypotheses of Theorem \ref{thm1} or Theorem \ref{thm2} is necessarily quasi-isometric.  Thus it makes sense to ask whether these results hold for general quasi-isometric laminations.\\

\noindent\textit{Quasi-isometric laminations with no compact leaves.}
Cannon and Thurston \cite{CanThurston} proved that the stable and unstable laminations of the suspension of a pseudo-Anosov homeomorphism of a closed surface are quasi-isometric, and each leaf is a plane or annulus in this case. In addition to these examples, Fenley \cite{Fenley2}
produced infinitely many examples of closed hyperbolic 3-manifolds with quasi-isometric laminations in which each leaf is an annulus, a mobius band, or a plane.  Note that Theorem \ref{thm2} implies that the examples of Cannon-Thurston and Fenley cannot have principal curvatures everywhere
in the interval $(-\delta_0 , \delta_0 )$.\\

One can also ask if we need to require that no complementary region is an interval bundle over a surface.\\

\noindent\textit{Small curvature laminations with simply-connected boundary leaves.}
Let $S$ be a closed totally geodesic embedded surface in a closed hyperbolic 3-manifold $M$.  Let $N(S) = S \times [0,1]$ be a closed embedded neighborhood of $S$ in $M$.  If the neighborhood $N(S)$ is small then the surfaces $S \times {t}$ will have small principal curvatures.
Since $\pi_1(S)$ is left-orderable, there exist faithful representations $\rho :\pi_1(S) \rightarrow Homeo([0,1])$ such that some points have trivial stabilizers (see Calegari \cite{Calegari}) The foliated bundle whose holonomy is $\rho$ has a leaf which is simply-connected.
Replace $N(S)$ with this foliated bundle.  We can blow up the simply-connected leaf and remove the interior to get a lamination which is $C^{\infty}$ close to the original (so that the leaves have small principal curvatures) and such that some boundary leaf is simply-connected.
See Calegari \cite{Calegari2} to see why the foliated bundle can be embedded in $M$ so that the leaves are smooth.  Note that this lamination has a complementary region which is an interval bundle over a surface.\\

\noindent\textit{Small curvature laminations with no compact leaves.}
One may also construct small curvature laminations in closed hyperbolic 3-manifolds with no compact leaves.  The author would like to thank Chris Leininger for describing the following construction.  The idea is to construct a small curvature branched surface in a closed hyperbolic 3-manifold which has an irrational point in the space of projective classes of measured laminations carried by the branched surface. A lamination corresponding to this irrational point will contain no compact leaves.  There are totally geodesic immersed closed surfaces in the figure-eight knot complement $M_8$ arbitrarily close to any plane in the tangent bundle (see Reid \cite{reid}).  Using this and the fact that $\pi_1(M_8)$ is LERF, one can find two such surfaces which lift to embedded surfaces $S_1$ and $S_2$ in a finite cover $M$ of $M_8$ which intersect in a non-separating (in both surfaces) simple closed geodesic $l$ at an arbitrarily small angle.  Flatten out the intersection to get a branched surface with small principal curvatures in which $S_1$ connects one side of $S_2$ to the other side.  The branched surface has three branch sectors (an annulus, $S_1 \setminus l$, and $S_2 \setminus l$) and one branch equation ($x_1 = x_2 +x_3$).  A solution to the branch equation in which two coordinates are not rationally related (e.g., $x_1 = 1/2$, $x_2 = 1/\pi$, $x_3 = 1/2 -1/\pi$)  will correspond to a lamination with no compact leaves which can be isotoped to have small principal curvatures.  Since the leaves do not have any cusps, we can fill the cusps of $M$ to get a small curvature lamination in a closed hyperbolic $3$-manifold with no compact leaves.



\section{Proof of Theorem \ref{thm1}}

Let $\epsilon >0$ be so small that if $P_1$, $P_2$, $P_3$ are three disjoint smoothly embedded planes in hyperbolic 3-space with principal curvatures in $(-1 , 1 )$ which intersect the same $\epsilon$-ball, then one of the $P_i$ separates the other two.

Let $\mathcal{L}$ be a codimension-one lamination in a finite volume hyperbolic 3-manifold $M$ such that the principal curvatures of each leaf are everywhere in the interval $(-\delta , \delta )$ for some $\delta \in (0,1)$.  Assume that no complementary region of $\mathcal{L}$ is an interval bundle over a surface.  Let $\tilde{\mathcal{L}}$ be the lift of $\mathcal{L}$ to $\h^3$.  Since every leaf of $\mathcal{L}$ has principal curvatures everywhere in $(-\delta ,\delta )$, the lamination $\mathcal{L}$ is a quasi-isometric lamination, and cannot be a foliation of $M$ by Fenley \cite{Fenley}.

Let $L_0$ be a boundary leaf of $\mathcal{L}$.  Suppose, for contradiction, that $\pi_1(L_0)$ is trivial, which implies that $L_0$ has infinite area.  Since $M$ is closed, $L_0$ must intersect some fixed compact ball in $M$ infinitely many times.  Thus given any integer $k$, we can find a point $y_k$ in $L_0$ such that the next leaf over on the boundary side of $L_0$ is within $1/k$ of $y_k$.

Let $\tilde{L}_0$ be a lift of $L_0$ to $\h^3$. Lift the points $y_k$ to a fixed fundamental domain of $\tilde{L}_0$ and call them $y_k$. Let $\tilde{L}_k$ be the next leaf over from $\tilde{L}_0$ which is within $1/k$ of $y_k$.  We now have a sequence of leaves $\tilde{L}_k$ in $\tilde{\mathcal{L}}$ on the boundary side of $\tilde{L}_0$ such that for each $k$ the distance from $\tilde{L}_k$ to $y_k$ is less than $1/k$, and there is no leaf of $\mathcal{L}$ between $\tilde{L}_0$ and $\tilde{L}_k$.  We also have that $\bd \tilde{L}_0 \neq \bd \tilde{L}_k$ for all $k$, because otherwise the region between $L_0$ and $L_k$ would be an interval bundle in the complement of $\mathcal{L}$.

Let $k$ be so large that $1/k < \epsilon /8$.  Since $\tilde{L}_k$ eventually diverges from $\tilde{L}_0$ we can find a point $x_k \in \tilde{L}_0$ such that the distance from $x_k$ to $\tilde{L}_k$ is exactly $\epsilon /8$.  Let $b_k$ be the $(\epsilon /32)$-ball tangent to $\tilde{L}_0$ at $x_k$ on the boundary side of $\tilde{L}_0$.

We will show that infinitely many of the balls $b_k$ are disjointly embedded in $M$, contradicting the fact that $M$ has finite volume.
Suppose that $\gamma (b_l) \cap b_k \neq \emptyset$ for some integers $l,k$ and some $\gamma$ in $\pi_1 (M)$. Note that $\gamma(\tilde{L}_0 ) \neq \tilde{L}_0$, since $L_0$ has trivial fundamental group.  Now $\tilde{L}_0$, $\tilde{L}_k$, and $\gamma(\tilde{L}_0 )$ all intersect some $\epsilon$-ball, so we must have that one of them separates the other two.  Since there are no leaves of $\tilde{\mathcal{L}}$ between $\tilde{L}_0$ and $\tilde{L}_k$, and $\gamma(\tilde{L}_0 )$ is closer to $x_k$ than $\tilde{L}_k$, we must have that $\tilde{L}_0$ separates $\tilde{L}_k$ and $\gamma(\tilde{L}_0 )$ (See figure 1(a) ).  Also note that $\tilde{L}_0$, $\tilde{L}_k$, and $\gamma(\tilde{L}_l )$ are all on the boundary side of $\gamma(\tilde{L}_0 )$ (i.e, the side which contains the ball $\gamma(b_l )$ ).

Now we will show no matter where $\gamma$ sends $\tilde{L}_l$, we get a contradiction.
We cannot have $\gamma(\tilde{L}_l ) = \tilde{L}_k$, because this would imply that $\gamma^{-1} (\tilde{L}_0 )$ separates $\tilde{L}_l$ and $\tilde{L}_0$.  Thus we have $\gamma(\tilde{L}_l ) \neq \tilde{L}_k$.

Since $\tilde{L}_0$, $\tilde{L}_k$, and $\gamma(\tilde{L}_l )$ all intersect some fixed $\epsilon$-ball, we must have that one of them separates the other two.
We cannot have that $\gamma(\tilde{L}_l )$ separates $\tilde{L}_0$ and $\tilde{L}_k$, because there are no leaves of $\tilde{\mathcal{L}}$ between $\tilde{L}_0$ and $\tilde{L}_k$ (See figure 1(b) ).
If $\tilde{L}_0$ separates $\tilde{L}_k$ and $\gamma(\tilde{L}_l )$, then $\gamma(\tilde{L}_l )$ is between $\tilde{L}_0$ and $\gamma(\tilde{L}_0 )$, so that $d(x_l , \tilde{L}_l ) = d(\gamma(x_l ), \gamma(\tilde{L}_l ) ) \le \epsilon /16$ which is a contradiction (See figure 1(c) ).  Thus $\tilde{L}_0$ cannot separate $\tilde{L}_k$ and $\gamma(\tilde{L}_l )$.
If $\tilde{L}_k$ separates $\tilde{L}_0$ and $\gamma(\tilde{L}_l )$, then $\gamma^{-1} (\tilde{L}_k )$ separates $\tilde{L}_0$ and $\tilde{L}_l$ which is a contradiction (See figure 1(d) ). Thus $\tilde{L}_k$ cannot separate $\tilde{L}_0$ and $\gamma(\tilde{L}_l )$.  We have shown that $\tilde{L}_l$ has nowhere to go under the map $\gamma$, so that $\gamma(b_l ) \cap \gamma(b_k ) = \emptyset$ for any integers $l,k$ and any $\gamma \in \pi_1 (M)$.
This implies that $M$ contains infinitely many disjoint $(\epsilon /32)$-balls, contradicting the fact that $M$ has finite volume. \hfill $\Box$\\

\begin{figure}[h]
\includegraphics[width=\textwidth]{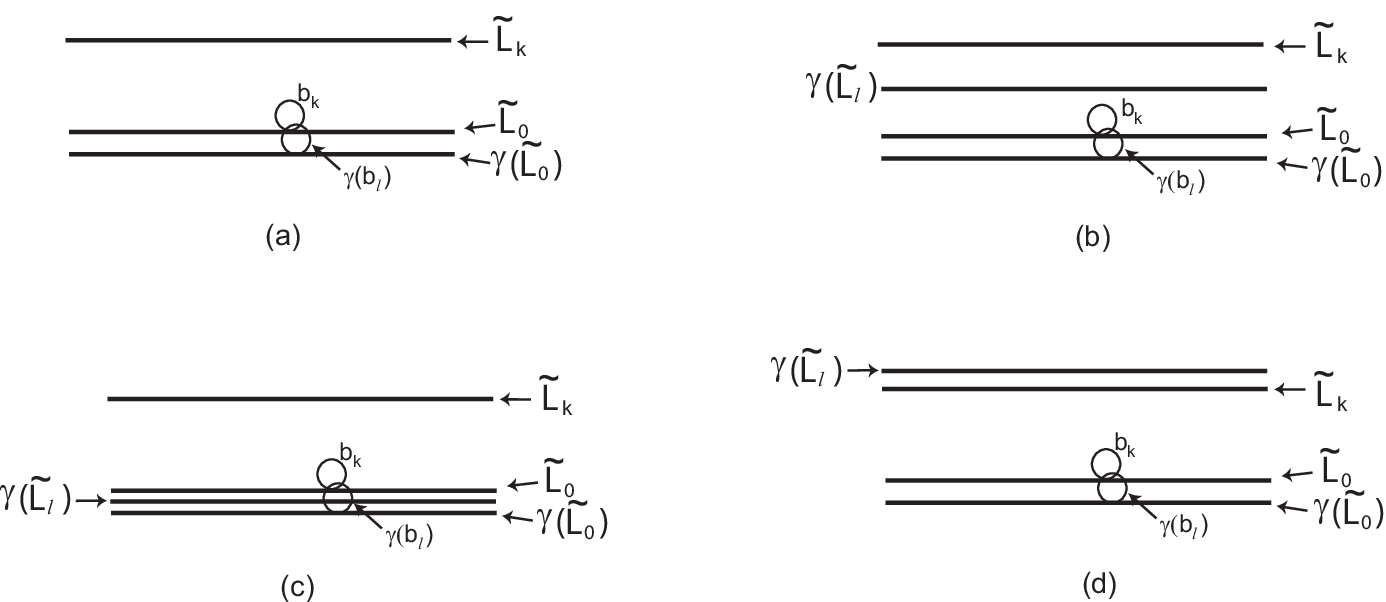}
\caption{(a) $\tilde{L}_0$ separates $\tilde{L}_k$ and $\gamma(\tilde{L}_0 )$.  (b) $\gamma(\tilde{L}_l )$ cannot separate $\tilde{L}_0$ and $\tilde{L}_k$.  (c) $\tilde{L}_0$ cannot separate $\tilde{L}_k$ and $\gamma(\tilde{L}_l )$. (d) $\tilde{L}_k$ cannot separate $\tilde{L}_0$ and $\gamma(\tilde{L}_l )$.}
\end{figure}

\section{Proof of Theorem \ref{thm2}}

Let $\epsilon > 0$ be so small that if $P_1$, $P_2$, $P_3$ are three disjoint smoothly embedded planes in hyperbolic 3-space with principal curvatures in $(-1,1)$ which intersect the same $\epsilon$-ball, then one of the $P_i$ separates the other two.  Let $\delta_0 > 0$ be so small that if a smooth curve $\gamma : (-\infty,\infty)\rightarrow \h^3$ in $\h^3$ with endpoints in $\bd\h^3$ has curvature at most $\delta_0$ at each point, then $\gamma(t)$ is in the $(\epsilon /2)$-neighborhood of the geodesic of $\h^3$ with the same endpoints.

Let $\mathcal{L}$ be a codimension-one lamination in a finite volume hyperbolic 3-manifold $M$ such that the principal curvatures of each leaf are everywhere in the interval $(-\delta_0 , \delta_0 )$.  Assume that no complementary region of $\mathcal{L}$ is an interval bundle over a surface.  Let $\tilde{\mathcal{L}}$ be the lift of $\mathcal{L}$ to $\h^3$.  As in the proof of Theorem \ref{thm1}, $\mathcal{L}$ cannot be a foliation. Let $L_0$ be a boundary leaf of $\mathcal{L}$.  Suppose, for contradiction, that $\pi_1(L_0)$ is cyclic, which implies that $L_0$ has infinite area.  Since $M$ is closed, $L_0$ must intersect some fixed compact ball in $M$ infinitely many times.  Also, by Theorem \ref{thm1}, we know that $\pi_1 (L_0)$ is nontrivial, so that $\pi_1 (L_0 ) \approx \Z$.

Let $\tilde{L}_0$ be a lift of $L_0$ to $\h^3$.  Since $L_0$ intersects a fixed compact ball in $M$ infinitely many times, we can find a sequence of points $y_k$ in $\tilde{L}_0$ such that the closest leaf of $\tilde{\mathcal{L}}$ to $y_k$ on the boundary side of $\tilde{L}_0$ is within $1/k$ of $y_k$.  Let $\tilde{L}_k$ be the leaf which is closest to $y_k$ on the boundary side of $\tilde{L}_0$.  Note that there is no leaf of $\tilde{\mathcal{L}}$ between $\tilde{L}_0$ and $\tilde{L}_k$.  We have $\bd \tilde{L}_0 \neq \bd \tilde{L}_k$ for all $k$, because the complement of $\mathcal{L}$ contains no interval bundle components.  We may assume that all $y_k$ are contained in a fixed fundamental domain $\mathcal{D}$ of $\tilde{L}_0$, and that $y_k$ converge to a point $y_{\infty} \in \bd\tilde{L}_0$.

For $k$ large enough we have $\bd \tilde{L}_0 \neq \bd \tilde{L}_k$ and $d(y_k , \tilde{L}_k ) \le \epsilon/8$, so that  we can find a point $x_k$ such that $d(x_k , \tilde{L}_k ) = \epsilon/8$.\\

\noindent\textbf{Case 1:}  We can choose the sequence of points $x_k \in \tilde{L}_0$ to be contained in a fixed fundamental domain $D$ of $\tilde{L}_0$ such that $x_k$ exit an end of $D$ whose projection to $M$ has infinite area.\\

Let $b_k$ be the $(\epsilon /32)$-ball tangent to $\tilde{L}_0$ at $x_k$ on the boundary side of $\tilde{L_0}$.  For $k$ large enough, say all $k$, the generator of $stab_{\pi_1 (M)} (\tilde{L}_0 )$ moves the center of $b_k$ a distance of at least $\epsilon$.  Thus we can assume that $\gamma(b_l ) \cap b_k = \emptyset$ for any integers $l,k$ and any $\gamma \in stab_{\pi_1 (M)} (\tilde{L}_0 )$.

We may now proceed as in the proof of Theorem \ref{thm1} to show that $\gamma(b_l ) \cap b_k = \emptyset$ for any integers $l,k$ and any $\gamma \in \pi_1 (M)$.  This again contradicts the fact that $M$ has finite volume.\\

\noindent\textbf{Case 2:}  We cannot choose the sequence of points $x_k$ as in Case 1.\\

If infinitely many of the leaves $\tilde{L}_k$ were distinct, then we would be able to find a sequence of points as described in Case 1.  Thus $\tilde{L}_k = \tilde{L}_+$ for some fixed leaf $\tilde{L}_+ \in \tilde{\mathcal{L}}$.

Let $U$ be the component of the complement in $\bd\tilde{L}_0$ of the fixed point(s) of the generator of $stab_{\pi_1 (M)} (\tilde{L}_0 )$ which contains the point $y_{\infty}$.  We will now show that $\bd\tilde{L}_+ \cap \bd\tilde{L}_0$ must contain $U$.

Suppose that $\bd\tilde{L}_+ \cap \bd\tilde{L}_0$ does not contain $U$. Since $d(y_k , \tilde{L}_+ ) < 1/k$ and $\bd\tilde{L}_+ \cap \bd\tilde{L}_0$ does not contain $U$, we can find a sequence of points $x_k$ in $\tilde{L}_0$ which converge to a point $x_{\infty} \in U$ with $d(x_k , \tilde{L}_+ ) = \epsilon /8$.  Since the point $x_{\infty}$ cannot be a fixed point of the generator of $stab_{\pi_1 (M)} (\tilde{L}_0 )$, a tail of the sequence $x_k$ must be contained in a fixed fundamental domain of $\tilde{L_0}$.  This contradicts the fact that we are in Case 2.  Thus $\bd\tilde{L}_+ \cap \bd\tilde{L}_0$ must contain $U$, hence must contain the fixed point(s) of the generator of $stab_{\pi_1 (M)} (\tilde{L}_0 )$.

If the generator of $stab_{\pi_1 (M)} (\tilde{L}_0 )$ is parabolic, then it has only one fixed point.  This implies that $\bd\tilde{L}_+ = \bd\tilde{L}_0$, giving us a contradiction.

If the generator of $stab_{\pi_1 (M)} (\tilde{L}_0 )$ is loxodromic, then we can argue as above to find a leaf $\bd\tilde{L}_-$ of $\tilde{\mathcal{L}}$ which contains the other component of complement in $\bd\tilde{L}_0$ of the fixed points of the generator of $stab_{\pi_1 (M)} (\tilde{L}_0 )$.
So $\bd\tilde{L}_+$ and $\bd\tilde{L}_-$ both contain the endpoints of the axis of the generator of $stab_{\pi_1 (M)} (\tilde{L}_0 )$. Since the principal curvatures of $\tilde{L}_0$, $\tilde{L}_+$, and $\tilde{L}_-$ are all in the interval $(-\delta_0 ,\delta_0 )$, and $\bd\tilde{L}_0$, $\bd\tilde{L}_+$, $\bd\tilde{L}_-$ all contain the endpoints of the axis of the generator of $stab_{\pi_1 (M)}$, we must have that $\tilde{L}_0$, $\tilde{L}_+$, and $\tilde{L}_-$ all intersect some fixed $\epsilon$-ball.  Thus one of the three separates the other two.  This gives us a contradiction since $\tilde{L}_+$ and $\tilde{L}_-$ are on the same side of $\tilde{L}_0$ (i.e, the boundary side) and there are no leaves of $\mathcal{L}$ between $\tilde{L}_0$ and $\tilde{L}_+$ or between $\tilde{L}_0$ and $\tilde{L}_-$. \hfill $\Box$\\

\textbf{Acknowledgements.}  This work was partially supported by the NSF grants DMS-0135345 and DMS-0602191.

%
%
%
\bibliographystyle{amsalpha}
\bibliography{tri}

\end{document}